\documentclass[12pt,reqno]{amsart}
\numberwithin{equation}{section}
\usepackage{amsmath}
\usepackage{amssymb}
\usepackage{amsthm,amsxtra}
\usepackage{latexsym}
\usepackage{cite}
\usepackage{enumerate}
\usepackage{enumitem}
\usepackage[mathscr]{eucal}
\usepackage[hyphens]{url}
\usepackage{mathtools}
\usepackage{multicol}
\usepackage[colorlinks=true]{hyperref}
\usepackage{cleveref}
\usepackage[top=2in, bottom=1.5in, left=1in, right=1in]{geometry}
\newtheorem{Def}{Definition}[section]
\newtheorem{Th}[Def]{Theorem}
\newtheorem{Ex}[Def]{Example}
\newtheorem{Lemma}[Def]{Lemma}
\newtheorem{Prop}[Def]{Proposition}
\newtheorem{Cor}[Def]{Corollary}
\newtheorem{Rem}[Def]{Remark}

{}

\crefname{Th}{Theorem}{Theorems}
\crefname{Lemma}{Lemma}{Lemmas}
\crefname{Cor}{Corollary}{Corollaries}
\crefname{Def}{Definition}{Definitions}
\crefname{Ex}{Example}{Examples}
\crefname{Prop}{Proposition}{Propositions}
\crefname{Prob}{Problem}{Problems}
\crefname{Rem}{Remark}{Remarks}
\crefname{Obs}{Observation}{Observations}
\crefname{Note}{Note}{Notes}
\DeclareMathOperator{\aut}{Sym(\omega)}

\DeclareMathOperator{\ks}{\kappa_{\sigma}}
\DeclareMathOperator{\gi}{G_\infty}
\DeclareMathOperator{\giz}{G_\infty^0}
\DeclareMathOperator{\gin}{G_\infty^{\neg 0}}

\begin{document}
\title[On L\'{e}vy type subgroups of the symmetric group of $\omega$]{On L\'{e}vy type subgroups of the symmetric group of $\omega$ associated with simple density functions}

\author[ N. Alam, P. Das]{Nur Alam$^\dag$, Pratulananda Das$^*$}
\address{\llap{\dag\,}Department of Mathematics, Mirza Ghalib College, Gaya, Bihar-823001, India}
\email{nurrejwan@gmail.com}
\address{\llap{*\,}Department of Mathematics, Jadavpur University, Kolkata-700032, West Bengal, India}
\email{pratulananda@yahoo.co.in}
%
%\address{\llap{\dag\,}Department of Mathematics, University of Gour
%Banga, Malda-732103, West Bengal, India}
%\email{debrajchandra1986@gmail.com}
%

%\thanks{The first author (as PI) is thankful to SERB(DST) for the CRG project (No. CRG/2022/000264) during the tenure of which this work has been done.}

\maketitle

\begin{abstract}
The L\'{e}vy group, introduced by L\'{e}vy in (L\'{e}vy (1951)), is defined by \[\mathcal{L} = \left\{\sigma\in \aut : \lim_n \frac{|\{k\in n : \sigma(k)\notin n\}|}{n} = 0\right\},\] where $\aut$ denotes the group of all permutations of $\omega$. This group happens to be closely connected with the notion of asymptotic density as it consists of permutations that map sets with an asymptotic density back to sets with the same density. Our aim is to replicate the construction by replacing $n$ in the denominator by $g(n)$ for a function $g:\omega\to [0,\infty)$ with $g(n)\to \infty$ (weight),  which leads us a class of L\'{e}vy type groups. Specifically we define simple  L\'{e}vy group with weight $g$ as \[\mathcal{L}_g = \left\{\sigma\in \aut : \lim_n \frac{|\{k\in n : \sigma(k)\notin n\}|}{g(n)} = 0\right\}.\]  As a natural consequence we first study its relationship with the notion of simple  density functions (which had been introduced in \cite{gd2} as a general class of density functions  by similar manner) paving way for the class of simple density functions and ideals $\mathcal{I}_g$. We have  $\mathcal{L}_g = \aut$ if and only if $n/g(n)\to 0$ (just like $\mathcal{I}_g = 2^{\omega}$ if and only if $n/g(n)\to 0$). The relation between simple  density functions and simple L\'{e}vy groups is in expected line. What seems more interesting is the comparative study of these groups. We start with the relation between $\mathcal{L}$ and $\mathcal{L}_g$ and constructing  suitable examples it is shown that strict inclusions and incomparability may occur. Finally we come up with an interesting observation that there are continuum many pairwise incomparable simple L\'{e}vy groups, each of which is incomparable with $\mathcal{L}$. 
\end{abstract}


\noindent{\bf\keywordsname{}:} {Asymptotic density, simple density, L\'{e}vy group, simple L\'{e}vy group, permutation group.}

\section{Introduction}
Let $\omega$ denote the set $\{0, 1, \dotsc\}$, and for $n\in \omega$, we use the identifications $n = \{0,1, \dotsc, n-1\}$ and $n+1 = \{0,1, \dotsc, n\}$. For a set $A\subseteq \omega$, its asymptotic density (or, density) is defined, whenever the limit exists, by \[d(A) = \lim_n \frac{|A\cap n|}{n}.\] Asymptotic density is one of the most natural ways of measuring the size of a subset of $\omega$. A classical problem in this area is to determine which permutations of $\omega$ preserve asymptotic density.

Let $\aut$ denote the group of all permutations of $\omega$, i.e. the symmetric group of $\omega$ (recall that for any set $X$, the symmetric group of $X$ is the group of all bijections from $X$ to $X$ \cite{cox19}). Historically there has been lots of interest in studying subgroups of $\aut$ (one can see \cite{LS, MN, Ne} for example).

One such instance is the study of L\'{e}vy group. For any $\sigma\in \aut$, we define $\ks(n) = |\{k\in n : \sigma(k)\notin n\}|$ for all $n\in \omega$. The set \[\mathcal{L} = \left\{\sigma\in \aut : \lim_n \frac{\ks(n)}{n} = 0\right\}\] forms a subgroup of $\aut$ and is called the L\'{e}vy group after T. Hida \cite{levy}. The L\'{e}vy group is closely connected with asymptotic density; see \cite{obata88, obata88-1, blum96, nath07, slez08}. Indeed, every permutation belonging to $\mathcal{L}$ preserves asymptotic density.

The definition of $\mathcal{L}$ is based on normalizing the boundary-crossing quantity $\ks(n)$ by the length $n$ of the initial segment. This naturally suggests replacing $n$ by a more general function and studying the resulting classes of permutations. Accordingly, for a function $g:\omega\to [0,\infty)$ with $g(n)\to \infty$, which we call a weight function, we introduce the simple L\'{e}vy group \[\mathcal{L}_g = \left\{\sigma\in \aut : \lim_n \frac{\ks(n)}{g(n)} = 0\right\}\] of weight $g$. Recall from \cite{gd1,gd2} that for a set $A\subseteq \omega$, whenever the limit exists, \[d_g(A) = \lim_n \frac{|A\cap n|}{g(n)}\] is called the simple density of $A$ of weight $g$ (the nomenclature is taken from subsequent set-theoretic investigations done in (see \cite{ak-1} for example).

The purpose of this article is to study the structure of these L\'{e}vy like groups or simple L\'{e}vy groups $\mathcal{L}_g$, their relationships with the classical L\'{e}vy group, their interrelationships, and their connections with simple densities.

We begin with some useful basic observations concerning simple densities. We show that the simple L\'{e}vy group $\mathcal{L}_g$ coincides with $\aut$ if and only if $n/g(n)\to 0$. Accordingly, we mainly focus on weights $g$ satisfying $n/g(n)\not\to 0$. We also prove that the simple density of weight $g$ is invariant under the corresponding simple L\'{e}vy group. The converse does not hold in general; however, we identify a class of weights for which the reverse implication is also valid.

Our next objective is to compare $\mathcal{L}_g$ with the classical L\'{e}vy group $\mathcal{L}$. The asymptotic behavior of the quotients $n/g(n)$ and $g(n)/n$ plays a central role in this comparison. We prove that $\mathcal{L} \subseteq \mathcal{L}_g$ if and only if the sequence $(n/g(n))$ is bounded, while the boundedness of $(g(n)/n)$ guarantees that $\mathcal{L}_g \subseteq \mathcal{L}$. Consequently, if both sequences are bounded, then $\mathcal{L}_g = \mathcal{L}$; in particular, this equality holds whenever $n/g(n)\to a$ for some $a\in (0,\infty)$. These conditions, however, are not necessary for equality: we construct a weight $g$ for which $\mathcal{L}_g = \mathcal{L}$ although $(g(n)/n)$ is unbounded. We construct examples of weight functions $g$ for which $\mathcal{L}_g \subsetneqq \mathcal{L}$, or $\mathcal{L} \subsetneqq \mathcal{L}_g$. Moreover, we construct a weight function $g$ such that $\mathcal{L}$ and $\mathcal{L}_g$ are incomparable. %It is to be noted that as some observations follow from the results established in Section~4 and in view of that, the examples given for simple density ideals in \cite{gd2}, but on the other hand there are results which actually have results in \cite{gd2} as their implications. As it is not a one way traffic, we present all the details to give a flavor as how results and examples can be established without going to set theory and make the article readable for those who do not want to go into details of set theory. Furthermore, most examples given in this section are carefully constructed to avoid similarity with examples given in \cite{gd2}.

%In \cite{slez08}, the classical L\'{e}vy group was characterized in terms of statistical convergence: a $\sigma\in \aut$ belongs to $\mathcal{L}$ if and only if the sequence $(\sigma(n)/n)$ is statistically convergent to $1$. To formulate a natural weighted analogue, we introduce statistical convergence with respect to a weight $g$. Motivated by the classical characterization, we investigate whether $\mathcal{L}_g$ can be characterized by any of the following conditions: the weighted statistical convergence of $(\sigma(n)/n)$ to $1$, the weighted statistical convergence of $(\sigma(n)/g(n))$ to $1$, the ordinary statistical convergence of $(\sigma(n)/g(n))$ to $1$, or the ordinary convergence of $(\sigma(n)/g(n))$ to $1$. We show by examples that none of these four natural characterizations is valid in general; indeed, for arbitrary weights $g$, both implications may fail in each case. Nevertheless, under additional assumptions on $g$, we obtain a positive characterization of $\mathcal{L}_g$ in terms of the weighted statistical convergence of the sequence $(\sigma(n)/n)$ to $1$.

As there is a class of special subgroups of $\aut$ now, it seems natural to investigate inclusion relations among simple L\'{e}vy groups. We first show that eventual domination of one weight function by a constant multiple of another yields the corresponding inclusion between their corresponding simple L\'{e}vy groups. In particular, if $g_1$ and $g_2$ are eventually comparable up to positive multiplicative constants, then $\mathcal{L}_{g_1} = \mathcal{L}_{g_2}$. This condition is sufficient but not necessary, as we construct weights $g_1$ and $g_2$ for which $(g_1(n)/g_2(n))$ is unbounded while $\mathcal{L}_{g_1} = \mathcal{L}_{g_2}$. We conclude this part by constructing a family $G$ of weights of cardinality $\mathfrak{c}$\footnote{$\mathfrak{c}$ denotes the cardinality of the continuum.} such that $\mathcal{L}_g$ is incomparable with $\mathcal{L}$ for every $g\in G$, and $\mathcal{L}_f$ and $\mathcal{L}_g$ are incomparable whenever $f,g\in G$ are distinct.

%Finally, we strengthen this observation by showing that there are actually continuum many pairwise non-isomorphic simple L\'{e}vy groups. This is where we use some set theory and existing results on simple density ideals from \cite{ak-1} by establishing a surprising fact that in certain instances simple L\'{e}vy groups just act like algebraic counterparts of simple density ideals.

\section{Basic observations on simple densities and definitions}
Let $\gi$ denote the collection of all weight functions, which will from now on simply be called weights, i.e. the functions $g: \omega \to [0,\infty)$ with $g(n) \to \infty$, $\giz$ denote the collection of all $g\in \gi$ with $n/g(n) \to 0$, and $\gin$ denote the collection of all $g\in \gi$ with $n/g(n) \not \to 0$. From \cite{gd1,gd2} we now recall the following variant of asymptotic density (or, density). Let $A \subseteq \omega$. Then the upper asymptotic density of $A$ of weight $g\in \gi$ is defined by \[\overline{d_g}(A) = \limsup_n \frac{|A\cap n|}{g(n)}\] and the lower asymptotic density of $A$ of weight $g$ is defined by \[\underline{d_g}(A) = \liminf_n \frac{|A\cap n|}{g(n)}.\]

If $\overline{d_g}(A) = \underline{d_g}(A)$, then the common value is called the asymptotic density of $A$ of weight $g$. Since its inception, a lot of work has been done on the associated ideals and these density functions have been renamed as simple density functions (see for example \cite{ak-1}). If $d_g(A)$ denotes the simple density of $A$ with weight $g$, then \[d_g(A) = \lim_n \frac{|A\cap n|}{g(n)}.\] Let $\mathcal{D}_g$ be the collection of all subsets of $\omega$ which admit the simple density of weight $g$. Note that $d_g(\omega) = 0$ if and only if $g\in \giz$ (see \cite{gd2}). For  $g(n)=n$, the notions reduce to the classical upper asymptotic density $\overline{d}$, lower asymptotic density $\underline{d}$, and asymptotic density $d$, respectively. We denote by $\mathcal{D}$ the collection of all subsets of $\omega$ having asymptotic density. Now we recall some of the features of the collection $\mathcal{D}_g$, mostly without proofs (compare with results given in \cite{gd2})

%Further recall from \cite{gd2} that $\mathcal{Z}_g =\{A \subset \omega : d_g(A) = 0\}$ forms a non-trivial proper free ideal (A family $\mathcal{I} \subseteq 2^Y$ of subsets of a nonempty set $Y$ is said to be an ideal in $Y$ if $(i)$ $A, B\in \mathcal{I}$ implies $A\cup B \in \mathcal{I}$ and $(ii)$ $A\in \mathcal{I}$, $B\subseteq A$ imply $B\in \mathcal{I}$. $\mathcal{I}$ is called proper and non-trivial if $\mathcal{I} \neq \{\emptyset\}$ and $Y \notin \mathcal{I}$. Further an admissible ideal $\mathcal{I}$ of $Y$  satisfies $\{x\}\in \mathcal{I}$ for each $x\in Y$. Such ideals are also called free ideals). This ideal has been named as simple density ideal after the corresponding density function (see \cite{ak-3, ak-2, ak-1}) and we will continue the same terminology. $\mathcal{Z} =\{A \subset \omega : d(A) = 0\}$ denotes, as usual, the ideal of all sets with asymptotic density zero.

\begin{Lemma}
\label{thm12}
$\mathcal{D}_g = P(\omega)$ if and only if $g\in \giz$.
\end{Lemma}

\begin{Lemma}
\label{lemma6}
If $\mathcal{D}\subseteq \mathcal{D}_g$, then either $\mathcal{D} = \mathcal{D}_g$ or $\mathcal{D}_g = P(\omega)$.
\end{Lemma}

\begin{Prop}
\label{prop4}
Let $g\in \gi$. Then $\mathcal{D} = \mathcal{D}_g$ if and only if there exists $a>0$ such that $n/g(n)\to a$. Moreover, whenever $n/g(n)\to a > 0$, we have $d_g(A)=ad(A)$ for every $A\in \mathcal{D}$.
\end{Prop}

\begin{proof}
Assume first that $\mathcal{D} = \mathcal{D}_g$. Then $\omega\in \mathcal{D}_g$, and so $\lim_n \frac{n}{g(n)}$ exists; say $\lim_n \frac{n}{g(n)} = a\geq 0$. Since $\mathcal{D} = \mathcal{D}_g$, by Lemma~\ref{thm12}, $a> 0$.

Conversely, assume that $n/g(n)\to a$ for some $a> 0$. For every $A\subseteq \omega$ and every $n\in \omega$, \[\frac{|A\cap n|}{g(n)} = \frac{|A\cap n|}{n}\cdot \frac{n}{g(n)}.\] It follows that $\mathcal{D}\subseteq \mathcal{D}_g$ and for every $A\in \mathcal{D}$, $d_g(A) = ad(A)$. Since $n/g(n)\to a > 0$, by Lemmas~\labelcref{thm12,lemma6}, we get $\mathcal{D} = \mathcal{D}_g$.
\end{proof}

\begin{Prop}
\label{prop15}
Let $g\in \gi$. If $(g(n)/n)$ converges, then $\mathcal{D}_g\subsetneqq \mathcal{D}$ if and only if $g(n)/n\to 0$.
\end{Prop}

\begin{proof}
The proof is straightforward and so is omitted.
\end{proof}

\begin{Th}
\label{thm5}
Let $A\in \mathcal{D}_g$. Then for any $\lambda$, $0\leq \lambda \leq d_g(A)$, there exists $B\subseteq A$ such that $d_g(B) = \lambda$.
\end{Th}

\begin{proof}
If $d_g(A) = 0$, then $\lambda = 0$, and taking $B = \emptyset$ we get $d_g(B) = 0 = \lambda$. Now assume that $d_g(A)> 0$, and put $p = \lambda/d_g(A) \in [0,1]$. Choose $A = \{a_0< a_1< a_2< \cdots\}$. Define \[B=\{a_m\in A : \lfloor (m+1)p\rfloor - \lfloor mp\rfloor = 1, m\in \omega\}.\] For each $n\in \omega$, let $x_n = |A\cap n|$. Then the elements of $A\cap n$ are exactly $a_0, a_1, a_2, \dotsc, a_{x_n-1}$. By the definition of $B$, the number of selected elements for $B$ among these is \[|B\cap n| = \sum_{m=0}^{x_n-1}\bigl(\lfloor (m+1)p\rfloor - \lfloor mp\rfloor\bigr) = \lfloor p x_n\rfloor.\] It follows that $0\leq p x_n - |B\cap n|< 1$ for all $n\in \omega$, and so \[0\leq \frac{px_n}{g(n)} - \frac{|B\cap n|}{g(n)} \leq \frac{1}{g(n)}\] for all $n\in \omega$. Since $1/g(n)\to 0$ and $x_n/g(n)\to  d_g(A)$, \[\frac{|B\cap n|}{g(n)} \to p d_g(A) = \frac{\lambda}{d_g(A)} d_g(A) = \lambda.\] Hence $d_g(B) = \lambda$.
\end{proof}

%\begin{Th}[{\cite[Theorem 1]{obata88}}]
%\label{thm2}
%Let $A\in \mathcal{D}$. If $d(A\Delta \sigma(A)) = 0$ for all $\sigma\in \mathcal{L}$, then $d(A) = 0$ or $1$.
%\end{Th}

%\begin{Cor}
%\label{cor4}
%Let $g\in \gi$ be such that $n/g(n) \to a\in (0,\infty)$, and let $A\in \mathcal{D}$. If $d_g(A\Delta \sigma(A)) = 0$ for all $\sigma\in \mathcal{L}$, then $d_g(A) = 0$ or $a$.
%\end{Cor}
%\begin{proof}
%By Proposition~\ref{prop4}, $\mathcal{D} = \mathcal{D}_g$ and $d_g(A) = ad(A)$. It follows that $d(A\Delta \sigma(A)) = 0$ for all $\sigma\in \mathcal{L}$. By Theorem~\ref{thm2}, $d(A) = 0$ or $1$. This gives $d_g(A) = 0$ or $a$.
%\end{proof}

\begin{Prop}
\label{prop6}
\hfill
\begin{enumerate}%[wide=0pt,label={\upshape(\arabic*)},leftmargin=*,ref={\theProp(\arabic*)}]
  \item\label{prop601} Let $A\in \mathcal{D}_g$. Then $A^c\in \mathcal{D}_g$ if and only if the sequence $(n/g(n))$ converges. Moreover, if $n/g(n) \to a\in [0,\infty)$, then $d_g(A^c) = a - d_g(A)$.
  \item\label{prop602} If $A, B \in \mathcal{D}_g$, then $\overline{d_g}(A\cup B) \leq d_g(A) + d_g(B)$.
  \item\label{prop603} If $A, B \in \mathcal{D}_g$ and $A\cap B\in \mathcal{Z}_g$, then $A\cup B\in \mathcal{D}_g$ and $d_g(A\cup B) = d_g(A) + d_g(B)$.
\end{enumerate}
\end{Prop}

\begin{Prop}
\label{prop7}
Let $A, B\in \mathcal{D}_g$. Then the following assertions are equivalent.
\begin{multicols}{2}
\begin{enumerate}%[label={\upshape(\arabic*)},leftmargin=*]
  \item $A\cup B\in \mathcal{D}_g$
  \item $A\cap B\in \mathcal{D}_g$
  \item $A\setminus B\in \mathcal{D}_g$
  \item $B\setminus A\in \mathcal{D}_g$.
\end{enumerate}
\end{multicols}
If any of the four conditions holds, then $d_g(A\cup B)=d_g(A) + d_g(B) - d_g(A\cap B)$.
\end{Prop}

\begin{Prop}
\label{prop8}
Let $A,B \subseteq \omega$ be such that $A\subseteq B$.
\begin{enumerate}%[wide=0pt,label={\upshape(\arabic*)},leftmargin=*,ref={\theProp(\arabic*)}]
  \item\label{prop801} Then $\underline{d_g}(A) \leq \underline{d_g}(B)$.
  \item\label{prop802} If $\overline{d_g}(B)$ is finite, then $\overline{d_g}(A) \leq \overline{d_g}(B)$.
  \item\label{prop803} If $A,B\in \mathcal{D}_g$, then $d_g(A)\leq d_g(B)$.
\end{enumerate}
\end{Prop}

Taking motivation from \cite{gd2,gd1} we now define the following generalized version of the classical L\'{e}vy group following the footsteps of \cite{levy}.

\begin{Def}
\label{def1}
Let $g\in \gi$. Then the subgroup \[\mathcal{L}_g = \left\{ \sigma\in \aut : \lim_n \frac{\ks(n)}{g(n)} = 0 \right\}\] of $\aut$ is said to be the simple L\'{e}vy group of weight $g$ (or just L\'{e}vy group of weight $g$).
\end{Def}

\begin{Th}
\label{thm13}
$\mathcal{L}_g = \aut$ if and only if $g\in \giz$.
\end{Th}

\begin{proof}
Suppose that $\mathcal{L}_g = \aut$. Let $A = \{2n : n\in \omega \setminus 2\omega\}$. We enumerate $2\omega$ and $A$ as $2\omega = \{r_0 < r_1 < \cdots\}$ and $\omega \setminus A = \{s_0 < s_1 < \cdots\}$. Clearly, $s_k \leq r_k$ for all $k\in \omega$. We now define $\sigma \in \aut$ by
\[\sigma(n) = \begin{cases}
2n, & \mbox{if } n\in \omega \setminus 2\omega\\
s_k, & \mbox{if } n = r_k \mbox{ for some } k\in \omega.
\end{cases}\]
Observe that
\[\ks(n) = \begin{cases}
\frac{n-2}{4} \text{ or } \frac{n+2}{4}, & \mbox{if } n\in 2\omega \\
\frac{n-3}{4} \text{ or } \frac{n-1}{4}, & \mbox{otherwise}.
\end{cases}\]
It follows that $\ks(n)/n \to 1/4$. Since $\mathcal{L}_g = \aut$, we obtain \[\frac{n}{g(n)} = \frac{\ks(n)}{g(n)}\cdot \frac{n}{\ks(n)} \to 0,\] i.e. $g\in \giz$.

Conversely, suppose that $g\in \giz$. Let $\sigma \in \aut$. For each $n\in \omega$, $\ks(n)/g(n) \leq n/g(n)$. It follows that $\ks(n)/g(n)\to 0$. Thus $\sigma\in \mathcal{L}_g$, and hence $\mathcal{L}_g = \aut$.
\end{proof}

\section{Simple L\'{e}vy groups}
\subsection{Connection with simple density functions}
In view of the construction of simple L\'{e}vy groups in Section~2 and the definition of simple density functions mentioned there, one would naturally expect connection between them (for a given weight $g$) in line of the already existing results showing connection between the asymptotic density and L\'{e}vy group. Also it is to be noted that as for $g(n) = n$ we get back the classical notions, all the existing results become special cases of the results established in this section.
%Furthermore, obviously as a natural recourse, we follow the lines of proofs from \cite{obata88}, in particular but there are necessary modifications.

\begin{Lemma}
\label{lemma2}
Let $\sigma\in \aut$ and $A\subseteq \omega$. Then for each $n\in \omega$, \[\bigl||\sigma(A)\cap n| - |A\cap n|\bigr| \leq \ks(n).\]
\end{Lemma}

\begin{proof}
For each $n\in \omega$, let $E_n = \{k\in n : \sigma(k)\notin n\}$ and $F_n = \{k\notin n : \sigma(k)\in n\}$. Since $\sigma\in \aut$, $|F_n| = |E_n| = \ks(n)$. Now for each $n\in \omega$, $\sigma^{-1}[n] = (n\setminus E_n) \cup F_n$, and the union is disjoint. Thus for each $n\in \omega$, \[|\sigma(A)\cap n| = |A\cap \sigma^{-1}[n]| = |A\cap (n\setminus E_n)| + |A\cap F_n|.\] On the other hand, for each $n\in \omega$, $|A\cap n| = |A\cap (n\setminus E_n)| + |A\cap E_n|$. Subtracting these two equalities, we obtain $|\sigma(A)\cap n| - |A\cap n| = |A\cap F_n|-|A\cap E_n|$ for all $n\in \omega$. Hence for each $n\in \omega$, \[\bigl||\sigma(A)\cap n| - |A\cap n|\bigr| = \bigl||A\cap F_n| - |A\cap E_n|\bigr| \leq \ks(n),\] because $|A\cap F_n| \leq |F_n| = \ks(n)$ and $|A\cap E_n| \leq |E_n| = \ks(n)$.
\end{proof}

\begin{Prop}
\label{prop9}
Simple density of weight $g$ is invariant under $\mathcal{L}_g$.
\end{Prop}

\begin{proof}
Let $\sigma\in \mathcal{L}_g$. We show that for every $A\in \mathcal{D}_g$, the simple density of $\sigma(A)$ of weight $g$ exists and $d_g(\sigma(A)) = d_g(A)$. Pick $A\in \mathcal{D}_g$. By Lemma~\ref{lemma2}, we get \[\left|\frac{|\sigma(A)\cap n|}{g(n)} - \frac{|A\cap n|}{g(n)}\right| \leq \frac{\ks(n)}{g(n)} \to 0.\] It follows that \[\frac{|\sigma(A)\cap n|}{g(n)} \to d_g(A).\] Thus the simple density of $\sigma(A)$ of weight $g$ exists and $d_g(\sigma(A)) = d_g(A)$.
\end{proof}

It follows from the preceding result that every member of $\mathcal{L}_g$ preserves the simple density of weight $g$ (thus extending the corresponding observation from \cite{obata88} about asymptotic density). Recall from \cite[Remark~2.3]{obata88} that there exists a $\sigma \in \aut$ such that $\sigma$ preserves the density but $\sigma \notin \mathcal{L}$. One can have a more general view of that remark from the next result.

\begin{Th}
\label{thm6}
If $g\in \gi$ satisfies $g(n+\lfloor g(n) \rfloor)/g(n) \to 1$, then there exists a $\sigma \in \aut$ such that $\sigma$ preserves the simple density of weight $g$, but $\sigma \notin \mathcal{L}_g$.
\end{Th}

\begin{proof}
We recursively choose a strictly increasing sequence $(N_k)$ as follows. Since $g(n)\to\infty$, for every $M\in\omega$, the set $\{g(n) : n\geq M\}$ has a minimum.
%Indeed, fix $M\in \omega$. Then there exists a $R> M$ such that $g(n)> g(M)$ for every $n\geq R$, and hence $\{g(n); n\geq M\}$ has the minimum at some point of the finite set $\{M, M+1, \dotsc, R-1\}$. Choose one such point and call it $N$. Then \[g(N) = \min\{g(n); n\geq M\}.\]
One can then choose $N_0$ for which $g(N_0) = \min\{g(n) : n\geq M_0\}$ for a sufficiently large $M_0$ and $N_0\geq M_0$. Put \[\ell_0 = \lfloor g(N_0)\rfloor, \qquad q_0 = N_0+\ell_0,\] and \[m_0 = \left\lfloor g(q_0)/4 \right\rfloor, \qquad r_0 = q_0+\lfloor g(q_0)\rfloor.\] After choosing $N_k$ and defining $r_k$, choose $M_{k+1}>r_k$ and select $N_{k+1}\geq M_{k+1}$ such that \[g(N_{k+1}) = \min\{g(n) : n\geq M_{k+1}\}.\] Then $N_{k+1}> r_k > N_k$, and for every $k\in\omega$, $g(n)\geq g(N_k)$ whenever $n\geq N_k$. By the hypothesis, we get \[\frac{g(q_k)}{g(N_k)} = \frac{g(N_k+\lfloor g(N_k)\rfloor)}{g(N_k)} \to 1.\] Applying the hypothesis once again at the points $q_k$, we obtain \[\frac{g(r_k)}{g(q_k)} = \frac{g(q_k+\lfloor g(q_k)\rfloor)}{g(q_k)} \to 1.\] Hence $g(r_k)/g(N_k) \to 1$. Since \[\frac{m_k}{\ell_k} \sim \frac{g(q_k)/4}{g(N_k)} = \frac{1}{4} \frac{g(q_k)}{g(N_k)} \to \frac{1}{4},\] we have $m_k\leq \ell_k$ for all sufficiently large $k$. By discarding finitely many initial terms of the sequence $(N_k)$, we may thus assume that $1\leq m_k\leq \ell_k$ for all $k\in\omega$. For each $k\in \omega$, define $I_k = [q_k-m_k, q_k) \cap \omega$ and $J_k = [q_k, q_k+m_k) \cap \omega$. Then for each $k\in \omega$, $I_k\cup J_k \subseteq [N_k,r_k)\cap \omega$. Moreover, the intervals $I_k\cup J_k$ are pairwise disjoint because $N_{k+1} > r_k$.

Define $\sigma\in \aut$ by
\[\sigma(n) = \begin{cases}
n+m_k, & \mbox{if }n\in I_k \mbox{ for some } k\in \omega\\
n-m_k, & \mbox{if }n\in J_k \mbox{ for some } k\in \omega\\
n, & \mbox{otherwise}.
\end{cases}\]
We claim that $d_g(\sigma(A)) = d_g(A)$ for all $A\in \mathcal{D}_g$. Let $A\in \mathcal{D}_g$, and let $d_g(A) = a$. Then $|A\cap n| - ag(n) = o(g(n))$ for all $n\in \omega$. It consequently follows that for each $k\in \omega$,
\begin{align}
\nonumber |A\cap [N_k,r_k)| & = |A\cap r_k| - |A\cap N_k|\\
\nonumber & = a(g(r_k) - g(N_k)) + o(g(r_k)) + o(g(N_k))\\
& = o(g(N_k)).
\label{eqn6}
\end{align}

\noindent\textbf{Case 1:} First suppose that $n\notin \bigcup_{k\in\omega} (I_k\cup J_k)$.

Since $\sigma^{-1}[n] = n$, we have $|\sigma(A)\cap n| = |A\cap \sigma^{-1}[n]| = |A\cap n|$. Therefore, \[ \frac{|\sigma(A)\cap n|}{g(n)} = \frac{|A\cap n|}{g(n)} \to a \text{ as } n\to \infty \text{ along } \omega \setminus \bigcup_{k\in \omega} (I_k\cup J_k).\]

\noindent\textbf{Case 2:} Now suppose that $n\in I_k\cup J_k$ for some $k\in \omega$.

If $n\in I_k$, then $n = q_k - m_k + t$, where $t\in m_k$. Since $\sigma^{-1}[n] = (q_k-m_k) \cup [q_k, q_k + t)$, it follows that \[|\sigma(A)\cap n| = |A\cap \sigma^{-1}[n]| = |A\cap (q_k-m_k)| + |A\cap [q_k,q_k+t)|.\] On the other hand, $|A\cap n| = |A\cap (q_k-m_k)| + |A\cap [q_k-m_k, q_k-m_k+t)|$. Consequently,
\begin{align*}
\bigl||\sigma(A)\cap n| - |A\cap n|\bigr| & = \bigl||A\cap [q_k,q_k+t)| - |A\cap [q_k-m_k, q_k-m_k+t)|\bigr|\\
& \leq |A\cap (I_k\cup J_k)|.
\end{align*}
If $n\in J_k$, then $n = q_k+t$, where $t\in m_k$. Since $\sigma^{-1}[n] = (q_k-m_k+t) \cup J_k$, we get \[|\sigma(A)\cap n| = |A\cap \sigma^{-1}[n]| = |A\cap (q_k-m_k)| + |A\cap [q_k-m_k, q_k-m_k+t)| + |A\cap J_k|.\] Also $|A\cap n| = |A\cap (q_k-m_k)| + |A\cap I_k| + |A\cap [q_k, q_k+t)|$. Therefore,
\begin{align*}
\bigl||\sigma(A)\cap n| - |A\cap n|\bigr| & = \bigl||A\cap [q_k+t, q_k+m_k)| - |A\cap [q_k-m_k+t, q_k)|\bigr|\\
& \leq |A\cap (I_k\cup J_k)|.
\end{align*}
Thus $\bigl||\sigma(A)\cap n| - |A\cap n|\bigr| \leq |A\cap [N_k,r_k)|$ as $I_k\cup J_k\subseteq [N_k,r_k)\cap \omega$. Evidently for $n\geq N_k$, we have  $g(n)\geq g(N_k)$. Therefore, by~\eqref{eqn6}, \[\frac{\bigl||\sigma(A)\cap n| - |A\cap n|\bigr| }{g(n)} \leq \frac{|A\cap [N_k,r_k)|}{g(N_k)} \to 0 \text{ as } n\to \infty \text{ along } \bigcup_{k\in \omega} (I_k\cup J_k).\] Consequently, \[\frac{|\sigma(A)\cap n|}{g(n)} - \frac{|A\cap n|}{g(n)} \to 0 \text{ as } n\to \infty \text{ along } \bigcup_{k\in \omega} (I_k\cup J_k),\] \[\text{i.e. } \frac{|\sigma(A)\cap n|}{g(n)} \to a \text{ as } n\to \infty \text{ along } \bigcup_{k\in \omega} (I_k\cup J_k)\]

Combining both cases, we can conclude that $d_g(\sigma(A))$ exists and $d_g(\sigma(A)) = d_g(A)$.

Finally, note that $\sigma \notin \mathcal{L}_g$. Indeed, \[\frac{\ks(q_k)}{g(q_k)} = \frac{m_k}{g(q_k)} = \frac{\lfloor g(q_k)/4 \rfloor}{g(q_k)} \to \frac{1}{4}.\] Hence $\sigma \notin \mathcal{L}_g$.
\end{proof}

\begin{Cor}
\label{cor2}
Let $g_1, g_2\in \gi$ be given by
\begin{enumerate}%[wide=0pt,label={\upshape(\arabic*)},leftmargin=*]
  \item $g_1(n) = \sqrt{n+1}$ for all $n\in \omega$;
  \item $g_2(n) = \ln\ln(n+3)$ for all $n\in \omega$.
\end{enumerate}
Then there exists $\sigma_i \in \aut$ for $i = 1, 2$ such that $\sigma_i$ preserves the simple density of weight $g_i$, but $\sigma_i \notin \mathcal{L}_{g_i}$, where $i=1, 2$.
\end{Cor}

\begin{Rem}
\label{rem2}
The condition $g(n+\lfloor g(n)\rfloor)/g(n) \to 1$ is not necessary for the conclusion of Theorem~\ref{thm6}. Indeed, let $g\in \gi$ be given by $g(n) = n+2+(-1)^n$ for all $n\in\omega$. Then $g(n+\lfloor g(n)\rfloor)/g(n) \to 2$. Indeed, if $n$ is even, then $g(n)=n+3$ and $g(n+\lfloor g(n)\rfloor) = g(2n+3) = 2n+4$; whereas if $n$ is odd, then $g(n) = n+1$ and $g(n+\lfloor g(n)\rfloor) = g(2n+1) = 2n+2$. Therefore, in either case, $g(n+\lfloor g(n)\rfloor)/g(n) \to 2$, which shows that the ratio does not converge to $1$.

On the other hand, observe that $n/g(n)\to 1$. It then follows that $\mathcal{D}_g = \mathcal{D}$ and $d_g(A) = d(A)$ for all $A\in \mathcal{D}$. Also for every $n\in \omega$ and every $\sigma\in \aut$, \[\frac{\ks(n)}{g(n)} = \frac{\ks(n)}{n} \cdot \frac{n}{g(n)},\] and hence $\mathcal{L}_g = \mathcal{L}$. From \cite[Remark~2.3]{obata88}, we have the implication that there exists a $\sigma\in \aut$ such that $\sigma$ preserves the simple density of weight $g$, but $\sigma \notin \mathcal{L}_g$.
\end{Rem}

\begin{Th}
\label{thm7}
Let $g\in \gi$ be such that there exists a strictly increasing sequence $(N_n)$ for which the following conditions hold.
\begin{enumerate}%[wide=0pt,label={\upshape(\arabic*)},leftmargin=*]
\item $\omega \setminus \{N_n : n\in \omega\}$ is infinite and $n/g(n)\to 0$ as $n\to \infty$ along $\omega \setminus \{N_n : n\in \omega\}$;
\item $\frac{\sum_{i\in n} g(N_i)}{g(N_n)} \to 0$.
\end{enumerate}
Then a $\sigma\in \aut$ belongs to $\mathcal{L}_g$ if and only if $\sigma$ preserves the simple density of weight $g$. If, in addition, $\limsup_n \frac{N_n}{g(N_n)} > 0$, then $\mathcal{L}_g \subsetneqq \aut$.
\end{Th}

\begin{proof}
Put $E = \{N_n : n\in\omega\}$ and $h_n = g(N_n)$ for every $n\in \omega$. Since $n/g(n)\to 0$ as $n\to \infty$ along $\omega \setminus E$, for every $A\subseteq \omega$, $|A\cap n|/g(n) \to 0$ as $n\to \infty$ along $\omega \setminus E$. Subsequently it follows that \[A\in \mathcal{D}_g \text{ if and only if } |A\cap N_n|/h_n \to 0,\tag{$*$}\] because $\omega \setminus E$ is infinite. Moreover, $d_g(A) = 0$ for every $A\in \mathcal{D}_g$. Again since $n/g(n)\to 0$ as $n\to \infty$ along $\omega \setminus E$, for every $\sigma\in \aut$, $\ks(n)/g(n) \to 0$ as $n\to \infty$ along $\omega \setminus E$. We can then conclude that \[\sigma\in \mathcal{L}_g \text{ if and only if } \ks(N_n)/h_n \to 0.\tag{$*_1$}\]

We already know that if $\sigma\in \mathcal{L}_g$, then $\sigma$ preserves the simple density of weight $g$. We would establish the reverse inclusion. Let $\sigma \in \aut$ be such that $\sigma$ preserves the simple density of weight $g$. If possible, suppose that $\sigma \notin \mathcal{L}_g$. Then by $(*_1)$, there exist an $\varepsilon > 0$ and an infinite set $N\subseteq \omega$ such that $\ks(N_n) \geq \varepsilon h_n$ for all $n\in N$. For each $n\in N$, define $K_n = \{k\notin N_n : \sigma(k)\in N_n\}$. Since $\sigma\in \aut$, $|K_n| = \ks(N_n)$ for every $n\in N$. By discarding finitely many elements of $N$, we may assume that $\sum_{i\in n} h_i< h_n$ for all $n\in N$. %it follows from the given condition (2).

Now choose recursively pairwise disjoint sets $B_n\subseteq K_n$, $n\in N$, such that $|B_n| = \lfloor\varepsilon h_n/4\rfloor$. We have, at the $n$-th stage, \[\left|\bigcup_{i\in N\cap n} B_i\right| \leq \frac{\varepsilon}{4} \sum_{i\in n} h_i < \frac{\varepsilon h_n}{4},\] whereas $|K_n|\geq \varepsilon h_n$. Set $A = \bigcup_{n\in N} B_n$. For every $r\in \omega$, the sets $B_n$ with $n\geq r$ are contained in $[N_r,\infty)$. Thus \[|A\cap N_r| \leq \sum_{i\in r} |B_i| \leq \frac{\varepsilon}{4} \sum_{i\in r} h_i.\] Consequently, \[\frac{|A\cap N_r|}{h_r} \leq \frac{\varepsilon}{4} \cdot \frac{\sum_{i\in r} h_i}{h_r} \to 0.\] By $(*)$, we subsequently obtain that $A\in \mathcal{D}_g$ and $d_g(A) = 0$.

On the other hand, for every $n\in N$, $\sigma(B_n) \subseteq N_n$. Hence \[\frac{|\sigma(A)\cap N_n|}{h_n} \geq \frac{|B_n|}{h_n} = \frac{\lfloor \varepsilon h_n/4 \rfloor}{h_n} \to \frac{\varepsilon}{4} \text{ as } n\to \infty \text{ along } N.\] Again by $(*)$, $\sigma(A) \notin \mathcal{D}_g$, i.e. $\sigma$ does not preserve the simple density of $A$ of weight $g$. Which is a contradiction. Therefore, $\sigma \in \mathcal{L}_g$.

Now take the additional assumption that $\limsup_n \frac{N_n}{h_n} > 0$. Then there exist a $\delta > 0$ and an infinite set $M\subseteq \omega$ such that $N_n\geq 4\delta h_n$ for all $n\in M$. Passing on to a subsequence of $M$, if necessary, we may assume that the intervals defined below are pairwise disjoint. Put $m_n = \lfloor \delta h_n\rfloor$, and define $U_n = [N_n-m_n, N_n) \cap \omega$ and $V_n = [N_n, N_n+m_n) \cap \omega$ for all $n\in M$. Define $\tau\in \aut$ by
\[\tau(k) = \begin{cases}
k+m_n, & \mbox{if } k\in U_n \mbox{ for some } n\in M\\
k-m_n, & \mbox{if } k\in V_n \mbox{ for some } n\in M\\
k, & \mbox{otherwise}.
\end{cases}\]
Then $\kappa_\tau(N_n) = m_n$ for all $n\in M$, and so \[\frac{\kappa_\tau(N_n)}{h_n} = \frac{\lfloor \delta h_n\rfloor}{h_n} \to \delta  \text{ as } n\to \infty \text{ along } M.\] Again by $(*_1)$, $\tau \notin \mathcal{L}_g$. This shows that $\mathcal{L}_g \subsetneqq \aut$.
\end{proof}

\begin{Cor}
\label{cor8}
For every $n\in \omega$, put $h_n = 2^{2^n}$ and $N_n = h_n^2$. If $g\in \gi$ is given by
\[g(k) = \begin{cases}
h_n, & \mbox{if } k = N_n \text{ for some } n\in \omega\\
k^2, & \mbox{otherwise},
\end{cases}\]
then a $\sigma\in \aut$ belongs to $\mathcal{L}_g$ if and only if $\sigma$ preserves the simple density of weight $g$. Moreover, $\mathcal{L}_g \subsetneqq \aut$.
\end{Cor}

\begin{proof}
Clearly, $h_{n+1} = h_n^2 = N_n$ for every $n\in \omega$. Let $E =\{N_n : n\in \omega\}$. Then $\omega \setminus E$ must be infinite. Moreover, for every $n\in \omega \setminus E$, we have $n/g(n) = n/n^2 = 1/n$. Therefore, $n/g(n) \to 0$ as $n\to \infty$ along $\omega \setminus E$.

Also $\sum_{i\in n} h_i \leq (n-1) h_{n-1}$ for all sufficiently large $n\in \omega$. Since $h_n = h_{n-1}^2$, we obtain \[\frac{\sum_{i\in n} h_i}{h_n} \leq \frac{(n-1) h_{n-1}}{h_{n-1}^2} = \frac{n-1}{h_{n-1}} \to 0.\] Thus the hypotheses of Theorem~\ref{thm7} are satisfied, and hence one can conclude that $\sigma\in \aut$ belongs to $\mathcal{L}_g$ if and only if $\sigma$ preserves the simple density of weight $g$.

Finally, $N_n/h_n = h_n^2/h_n = h_n \to \infty$. In particular, $\limsup_n \frac{N_n}{h_n} > 0$. Applying the final assertion of Theorem~\ref{thm7}, we can observe that $\mathcal{L}_g \subsetneqq \aut$.
\end{proof}

The proof of the following result follows the same lines as that of \cite[Lemma~2.4]{obata88} with the necessary modifications.

\begin{Th}
\label{thm1}
Let $A = \{a_0< a_1< \cdots\}$ and $B = \{b_0< b_1< \cdots\} \in \mathcal{D}_g$ be such that $d_g(A) = d_g(B)$. Put $A^c = \{a_0^\prime< a_1^\prime< \cdots\}$ and $B^c = \{b_0^\prime< b_1^\prime< \cdots\}$. Assume that $A, A^c, B$ and $B^c$ are infinite. Let $\sigma\in \aut$ be given by $\sigma(a_n) = b_n$, $\sigma(a_n^\prime) = b_n^\prime$ for all $n\in \omega$. Then $\sigma\in \mathcal{L}_g$.
\end{Th}

\begin{proof}
For each $n\in \omega$, let $\alpha(n) = |A\cap n|$, $\beta(n) = |B\cap n|$, $\alpha^\prime(n) = |A^c\cap n|$ and $\beta^\prime(n) = |B^c\cap n|$. Clearly, $\alpha^\prime(n) - \beta^\prime(n) = \beta(n) - \alpha(n)$ for all $n\in \omega$. Since $d_g(A) = d_g(B)$, we have $\lim_n \frac{\alpha(n)-\beta(n)}{g(n)} = 0$.

Observe that for each $n\in \omega$,
\begin{equation}
\left.
\begin{array}{c}
a_{\alpha(n)}\leq n<a_{\alpha(n)+1},\;b_{\beta(n)}\leq n<b_{\beta(n)+1}\\
a_{\alpha(n)}^\prime\leq n<a_{\alpha(n)+1}^\prime,\;b_{\beta(n)}^\prime\leq n<b_{\beta(n)+1}^\prime
\end{array}\right\}
\label{eqn3}
\end{equation}

Now \[\{k\in n : \sigma(k) \notin n\} = \{a_k\in A\cap n : \sigma(a_k) \notin n, k\in \omega\} \cup \{a_k^\prime\in A^c\cap n : \sigma(a_k^\prime) \notin n, k\in \omega\}\] for all $n\in \omega$. Subsequently it follows that
\begin{eqnarray*}
\ks(n) &=& |\{a_k\in A\cap n : \sigma(a_k) \notin n, k\in \omega\}| + |\{a_k^\prime\in A^c\cap n : \sigma(a_k^\prime) \notin n, k\in \omega\}|\\
&=& |\{a_k\in A\cap n : b_k \notin n, k\in \omega\}| + |\{a_k^\prime \in A^c\cap n : b_k^\prime \notin n, k\in \omega\}|\\
&=& |\{a_k\in A\cap n : b_k\notin B\cap n, k\in \omega\}| + |\{a_k^\prime \in A^c\cap n : b_k^\prime \notin B^c\cap n, k\in \omega\}|\\
&=& \max\{0,\alpha(n) - \beta(n)\} + \max\{0, \alpha^\prime(n) - \beta^\prime(n)\} \; [\text{by \eqref{eqn3}}]\\
&=& |\alpha(n)-\beta(n)|
\end{eqnarray*} for all $n\in \omega$. This implies that
$\ks(n)/g(n)\to 0$ and so  $\sigma\in \mathcal{L}_g$.
\end{proof}

\begin{Cor}
\label{cor3}
Let $A, B\in \mathcal{D}_g$ be such that $d_g(A) = d_g(B)$. Assume that $A$, $A^c$, $B$ and $B^c$ are infinite. Then there exists a $\sigma\in \mathcal{L}_g$ such that $\sigma(A) = B$.
\end{Cor}

In \cite{obata88}, Obata had proved that the asymptotic density can uniquely be determined by finite additivity, invariance under the L\'{e}vy group $\mathcal{L}$, and normalization at $\omega$. More precisely:

\begin{Th}[{\cite[Theorem 2]{obata88}}]
\label{thm15}
If $\gamma: \mathcal{D}\to [0,\infty)$ satisfies
\begin{enumerate}%[wide=0pt,label={\upshape(\arabic*)},leftmargin=*]
\item (additivity) $\gamma(A\cup B) = \gamma(A)+\gamma(B)$ for every pair of disjoint sets $A,B\in \mathcal{D}$;
\item (invariance) $\gamma(\sigma(A)) = \gamma(A)$ for every $A\in \mathcal{D}$ and every $\sigma\in \mathcal{L}$;
\item (normalization) $\gamma(\omega) = 1$,
\end{enumerate}
then $\gamma = d$.
\end{Th}

Since the simple density $d_g$ does not satisfy the normalization property in general, it cannot be characterized in the same manner as the classical asymptotic density. If $d_g(\omega) = a\in (0,\infty)$, then we can obtain a characterization of the simple density as follows which is itself a consequence of the classical result.

\begin{Cor}
\label{cor9}
Assume that $d_g(\omega) = a \in (0,\infty)$. If $\gamma: \mathcal{D}_g\to [0,\infty)$ satisfies
\begin{enumerate}%[wide=0pt,label={\upshape(\arabic*)},leftmargin=*]
\item $\gamma(A\cup B)=\gamma(A)+\gamma(B)$ for every pair of disjoint
sets $A,B\in\mathcal{D}_g$;
\item $\gamma(\sigma(A))=\gamma(A)$ for every $A\in\mathcal{D}_g$ and
every $\sigma\in\mathcal{L}_g$;
\item $\gamma(\omega) = a$,
\end{enumerate}
then $\gamma = d_g$.
\end{Cor}

\begin{Rem}
\label{rem4}
If $d_g(\omega) = 1$, then the normalization condition in the preceding corollary becomes $\gamma(\omega) = 1$ and consequently we obtain the classical characterization theorem.
\end{Rem}

Let $S$ be a set and $\Sigma$ be a nonempty family of subsets of $S$. Recall that $\Sigma$ is said to be a field if it satisfies the following conditions.
\begin{enumerate}
  \item $\emptyset\in \Sigma$.
  \item If $A\in \Sigma$, then $A^c\in \Sigma$.
  \item If $\mathcal{A}$ is a finite subset of $\Sigma$, then $\bigcup \mathcal{A}\in \Sigma$.
\end{enumerate}
It is to be noted that if $(1)$ holds, then $(2)$ is equivalent to the assertion that if $\mathcal{A}$ is a finite subset of $\Sigma$, then $\bigcap \mathcal{A}\in \Sigma$.

Let $\widetilde{\mathcal{D}}_g$ denote the smallest field of subsets of $\omega$ containing $\mathcal{D}_g$. From Proposition~\ref{prop9}, we can conclude that $\mathcal{D}_g$ is stable under $\mathcal{L}_g$, i.e. $\sigma \mathcal{D}_g = \mathcal{D}_g$. We now show that the field generated by $\mathcal{D}_g$ is also stable under $\mathcal{L}_g$.

\begin{Prop}
\label{prop10}
The field $\widetilde{\mathcal{D}}_g$ is stable under $\mathcal{L}_g$.
\end{Prop}

\begin{proof}
The proof is straightforward and so is omitted.
%Fix $\sigma\in \mathcal{L}_g$, and define $\Sigma_\sigma = \{A\subseteq \omega : \sigma(A)\in \widetilde{\mathcal{Z}}_g\}$. Since $\widetilde{\mathcal{Z}}_g$ is a field and $\sigma\in \aut$, the family $\Sigma_\sigma$ is also a field. Since $\sigma\mathcal{Z}_g = \mathcal{Z}_g$, we have $\mathcal{Z}_g\subseteq \Sigma_\sigma$. Since $\widetilde{\mathcal{Z}}_g$ is the smallest field containing $\mathcal{Z}_g$, we obtain $\widetilde{\mathcal{Z}}_g \subseteq \Sigma_\sigma$. It follows that $\sigma\widetilde{\mathcal{Z}}_g \subseteq \widetilde{\mathcal{Z}}_g$, and hence $\sigma\widetilde{\mathcal{D}}_g =
%\widetilde{\mathcal{Z}}_g$.
\end{proof}

Weakening the assumptions of Theorem~\ref{thm15}, Obada \cite{obata88} had obtained the following characterization of asymptotic density.

\begin{Th}[{\cite[Theorem 3]{obata88}}]
\label{thm9}
Let $\gamma$ be a real-valued function defined on $\mathcal{D}$. If $\gamma$ can be extended to a real-valued function $\widetilde{\gamma}$ on the field $\widetilde{\mathcal{D}}$ (where $\widetilde{\mathcal{D}}$ denotes the smallest field of subsets of $\omega$ containing $\mathcal{D}$) such that \begin{enumerate}%[wide=0pt,label={\upshape(\arabic*)},leftmargin=*]
\item $\widetilde{\gamma}$ is finitely additive;
\item $\widetilde{\gamma}$ is of bounded variation;
\item $\widetilde{\gamma}$ is invariant under $\mathcal{L}$,
\end{enumerate}
then $\gamma = cd$, where $c = \gamma(\omega)$.
\end{Th}

Because of the lack of normalization, the above result also does not hold in general for the simple density $d_g$. If $d_g(\omega) = a\in(0,\infty)$, then, weakening the assumptions of Corollary~\ref{cor9}, we obtain the following characterization of the simple density of weight $g$ as a consequence of Theorem~\ref{thm9}.

\begin{Cor}
\label{cor10}
Assume that $d_g(\omega) = a\in (0,\infty)$. Let $\gamma$ be a real-valued function defined on $\mathcal{D}_g$. If $\gamma$ can be extended to a real-valued function $\widetilde{\gamma}$ on the field $\widetilde{\mathcal{D}}_g$ such that
\begin{enumerate}%[wide=0pt,label={\upshape(\arabic*)},leftmargin=*]
\item $\widetilde{\gamma}$ is finitely additive;
\item $\widetilde{\gamma}$ is of bounded variation;
\item $\widetilde{\gamma}$ is invariant under $\mathcal{L}_g$,
\end{enumerate}
then $\gamma = cd_g/a$, where $c = \gamma(\omega)$.
\end{Cor}

%\begin{Cor}
%\label{cor11}
%Assume that $d_g(\omega) = a\in (0,\infty)$. Let $\gamma$ be a real-valued function defined on $\mathcal{Z}_g$ such that $\gamma(\omega) = a$. If $\gamma$ can be extended to a real-valued function $\widetilde{\gamma}$ on the field $\widetilde{\mathcal{Z}}_g$ such that
%\begin{enumerate}%[wide=0pt,label={\upshape(\arabic*)},leftmargin=*]
%\item $\widetilde{\gamma}$ is finitely additive;
%\item $\widetilde{\gamma}$ is of bounded variation;
%\item $\widetilde{\gamma}$ is invariant under $\mathcal{L}_g$,
%\end{enumerate}
%then $\gamma = d_g$.
%\end{Cor}

\subsection{Relations between simple L\'{e}vy groups for different weight functions}
We start with the observation that unlike the classical case \cite{obata88} now we do not have only one subgroup of $\aut$ (i.e. the L\'{e}vy group which corresponds to asymptotic density) but actually a class of subgroups arising out of different weight functions (i.e. different simple density functions). This class obviously contains the L\'{e}vy group but we want to know whether it contains more and if so, how many.

In this section, our main focus is on the interrelationships between simple L\'{e}vy groups when different weight functions occur. In order to do that we can look at the interrelationships between simple density ideals (this will be taken up in Section~4) and this surprising connection enables us to conclude some observations about simple L\'{e}vy groups from the existing results for simple density ideals but not all. In particular the study of the relation with a simple L\'{e}vy group of weight $g$ with the L\'{e}vy group (when one is subset of the other and when they are incomparable etc.) does not follow from the relation of simple density ideals with asymptotic density zero ideal. Hence we start with several examples and results to understand the relation of simple L\'{e}vy groups with the classical L\'{e}vy group first, before proceeding on to investigate their interrelationships. Readers should note that in most of the examples given below we do not take help of the weight functions already used in articles like \cite{gd2} or \cite{ak-1}.

\begin{Th}
\label{thm11}
Let $g\in \gi$. Then $\mathcal{L} \subseteq \mathcal{L}_g$ if and only if the sequence $(n/g(n))$ is bounded.
\end{Th}

\begin{proof}
Let $\mathcal{L} \subseteq \mathcal{L}_g$. If possible, suppose that $(n/g(n))$ is unbounded. Then there exists a strictly increasing sequence $(n_k)$ in $\omega$ such that $n_k/g(n_k)\to \infty$. For each $k\in \omega$, put $\ell_k = \lfloor g(n_k)\rfloor$. Then $\ell_k\to \infty$ and $\ell_k/n_k \to 0$. Passing to a subsequence, if necessary, we may assume that $n_{k+1}> n_k + \ell_k + \ell_{k+1}$ for all $k\in \omega$, i.e. the intervals $I_k\cup J_k$ are pairwise disjoint, where $I_k = [n_k-\ell_k,n_k) \cap \omega$ and $J_k = [n_k,n_k+\ell_k) \cap \omega$. Define $\sigma\in \aut$ by
\[\sigma(n) = \begin{cases}
n+\ell_k, & \mbox{if } n\in I_k \mbox{ for some } k\in \omega\\
n-\ell_k, & \mbox{if } n\in J_k \mbox{ for some } k\in \omega\\
n, & \mbox{otherwise}.
\end{cases}\]
If $n\notin \bigcup_{k\in \omega} (I_k\cup J_k)$, then $\ks(n) = 0$. If $n\in I_k\cup J_k$ for some $k\in \omega$, then $\ks(n)\leq \ell_k$. Since $n\geq n_k-\ell_k$, we obtain \[\frac{\ks(n)}{n}\leq \frac{\ell_k}{n_k-\ell_k}.\] Thus $\ks(n)/g(n)\to 0$, i.e. $\sigma\in \mathcal{L}$.

On the other hand, $\ks(n_k) = \ell_k$ for all $k\in \omega$, and therefore \[\frac{\ks(n_k)}{g(n_k)} = \frac{\lfloor g(n_k)\rfloor}{g(n_k)} \to 1.\]
Thus $\sigma \notin \mathcal{L}_g$, contradicting $\mathcal{L} \subseteq \mathcal{L}_g$. Hence the sequence $(n/g(n))$ must be bounded.

The converse follows immediately from the equality \[\frac{\ks(n)}{g(n)} = \frac{\ks(n)}{n} \cdot \frac{n}{g(n)}\] for every $n\in \omega$.
\end{proof}

\begin{Prop}
\label{prop3}
Let $g\in \gi$. If the sequence $(g(n)/n)$ is bounded, then $\mathcal{L}_g\subseteq \mathcal{L}$.
\end{Prop}
\begin{proof}
Let $\sigma \in \mathcal{L}_g$. Now for each $n\in \omega$,
\[\frac{\ks(n)}{n} = \frac{\ks(n)}{g(n)} \cdot \frac{g(n)}{n}.\] This gives $\ks(n)/n\to 0$, i.e. $\sigma \in \mathcal{L}$. Thus $\mathcal{L}_g \subseteq \mathcal{L}$.
\end{proof}

\begin{Cor}
\label{cor1}
Let $g\in \gi$.
\begin{enumerate}%[wide=0pt,label={\upshape(\arabic*)},leftmargin=*]
  \item If both the sequences $(n/g(n))$ and $(g(n)/n)$ are bounded, then $\mathcal{L}_g = \mathcal{L}$.
  \item If $n/g(n) \to a\in (0,\infty)$, then $\mathcal{L}_g = \mathcal{L}$.
\end{enumerate}
\end{Cor}

\begin{Lemma}
\label{lemma5}
Let $\sigma\in \aut$. Then for every $m, n\in \omega$, $|\ks(m) - \ks(n)| \leq |m-n|$.
\end{Lemma}

\begin{proof}
We first show that $|\ks(n+1) - \ks(n)| \leq 1$ for all $n\in \omega$. For each $n\in \omega$, let $A_n = \{k\in n : \sigma(k)\notin n\}$ for all $n\in \omega$, and $B = \{k\in n : \sigma(k)\notin n+1\}$. Then for every $n\in \omega$, \[A_{n+1} = B\cup (\{n\}\cap \{k\in \omega : \sigma(k)\notin n+1\}),\] and so $|A_{n+1}| = |B| + b$, where $b\in \{0,1\}$.

On the other hand, for every $n\in \omega$, $A_n = B\cup \{k\in n : \sigma(k)= n\}$. Since $\sigma$ is a permutation, \[|\{k\in n : \sigma(k)= n\}|\in \{0,1\}.\] It follows that $|A_n| = |B| + c$, where $c\in\{0,1\}$. Thus \[\ks(n+1) - \ks(n) = b-c\in \{-1,0,1\}\] for all $n\in \omega$. Hence $|\ks(n+1) - \ks(n)|\le 1$ for all $n\in \omega$.

Now let $m, n\in \omega$. Without loss of generality, assume that $m< n$. Then \[|\ks(n) - \ks(m)| \leq \sum_{j=m}^{n-1} |\ks(j+1) - \ks(j)|.\] Using the above estimate, we get \[|\ks(n) - \ks(m)| \leq \sum_{j=m}^{n-1} 1 = n-m.\] Therefore, $|\ks(m) - \ks(n)| \leq |m-n|$.
\end{proof}

\begin{Ex}
\label{ex8}
There exists a $g\in \gin$ such that $\mathcal{L}_g = \mathcal{L}$, but the sequence $(g(n)/n)$ is unbounded.
\end{Ex}

\begin{proof}
For each $k\in \omega$, let $n_k = 2^{2^k}$. Define $g\in \gin$ by
\[g(n) = \begin{cases}
n_k^2, & \mbox{if } n = n_k \mbox{ for some } k\in \omega\\
n, & \mbox{otherwise}.
\end{cases}\]
Since $g(n_k)/n_k = n_k\to \infty$, $(g(n)/n)$ is unbounded. We claim that $\mathcal{L}_g = \mathcal{L}$. As $(n/g(n))$ is bounded, by Theorem~\ref{thm11}, $\mathcal{L} \subseteq \mathcal{L}_g$. Now let $\sigma\in \mathcal{L}_g$. If possible, suppose that $\sigma \notin \mathcal{L}$. Then there exists $\varepsilon > 0$ such that $A = \{n\in \omega : \ks(n)\geq \varepsilon n\}$ is infinite. Choose $B = A\cap (\omega \setminus \{n_k : k\in \omega\})$. Obviously, \[\frac{\ks(n)}{g(n)} = \frac{\ks(n)}{n} \geq \varepsilon\] for all $n\in B$. It follows that $B$ is finite as $\sigma\in \mathcal{L}_g$. Thus $A\cap \{n_k : k\in \omega\}$ is infinite, i.e. for infinitely many $k\in \omega$, we have $\ks(n_k)\geq \varepsilon n_k$. By Lemma~\ref{lemma5}, $\ks(n_k+1)\geq \ks(n_k)-1\geq \varepsilon n_k-1$ for all $k\in \omega$. Further as $g(n_k+1) = n_k+1$ for all $k\in \omega$, so \[\frac{\ks(n_k+1)}{g(n_k+1)} \geq \frac{\varepsilon n_k-1}{n_k+1} \to \varepsilon.\] This contradicts the fact that $\sigma\in \mathcal{L}_g$. Therefore, we must have $\sigma\in \mathcal{L}$, and so $\mathcal{L}_g \subseteq \mathcal{L}$. Thus $\mathcal{L}_g = \mathcal{L}$.
\end{proof}

\begin{Ex}
\label{ex3}
There exists $g\in \gin$ such that $\mathcal{L}_g \subsetneqq \mathcal{L}$.
\end{Ex}

\begin{proof}
Let $g\in \gin$ be defined by $g(n) = \sqrt{n}$ for all $n \in \omega$. Clearly, $g(n)/n \to 0$, and hence by Proposition~\ref{prop3}, $\mathcal{L}_g \subseteq \mathcal{L}$. It remains to show that the inclusion is proper. Let $A \subseteq \omega$ be the collection of all perfect squares and $B = \{n^2 : n\in A\}$. We enumerate $\omega \setminus A$ and $\omega \setminus B$ as $\omega \setminus A = \{p_0< p_1< \cdots\}$ and $\omega \setminus B = \{q_0 < q_1 < \cdots\}$. Clearly, $p_k> q_k$ for all sufficiently large $k\in \omega$. Let us now define $\sigma \in \aut$ by
\[\sigma(n) = \begin{cases}
n^2, & \mbox{if } n\in A\\
q_k, & \mbox{if } n = p_k \mbox{ for some } k\in \omega.
\end{cases}\]
Observe that
\begin{equation}
\ks(n) = \lfloor \sqrt{n-1}\rfloor - \left\lfloor \sqrt{\lfloor \sqrt{n-1} \rfloor}\right\rfloor
\label{eqn2}
\end{equation} for all sufficiently large $n\in \omega$. It then follows that \[\lim_n \frac{\ks(n)}{n} = \lim_n \frac{\lfloor \sqrt{n-1}\rfloor - \left\lfloor \sqrt{\lfloor \sqrt{n-1} \rfloor} \right\rfloor}{n} = 0.\] and so $\sigma \in \mathcal{L}$. From \eqref{eqn2}, we further have \[\lim_n \frac{\ks(n)}{g(n)} = \lim_n \frac{\lfloor\sqrt{n-1}\rfloor - \left\lfloor \sqrt{\lfloor \sqrt{n-1} \rfloor} \right\rfloor}{\sqrt{n}} = 1.\] Consequently, $\sigma \notin \mathcal{L}_g$. Thus $\mathcal{L}_g \subsetneqq \mathcal{L}$.
\end{proof}

\begin{Ex}
\label{ex03}
There exists $g\in \gin$ for which $\mathcal{L}\subsetneqq \mathcal{L}_g \subsetneqq \aut$.
\end{Ex}

\begin{proof}
Let $(x_n)$ be given by $x_n = 2^{2^n}$ for all $n\in \omega$. Define $y_0 = 0$, $y_{n+1} = y_n + 2 x_n$, and put $M_n = [y_n, y_n + x_n) \cap \omega$, $N_n = [y_n + x_n, y_{n+1}) \cap \omega$. Then $\{M_n, N_n : n\in \omega\}$ forms a partition of $\omega$. Now define $g\in \gin$ by
\[g(n)= \begin{cases}
n^2, & \mbox{if } n\in M_{2m}\cup N_{2m} \mbox{ for some } m\in \omega\\
n, & \mbox{if } n\in M_{2m+1}\cup N_{2m+1} \mbox{ for some } m\in \omega.
\end{cases}\]
Clearly, $(n/g(n))$ is bounded and hence by Theorem~\ref{thm13}, $\mathcal L_g \subsetneqq \aut$, whereas by Theorem~\ref{thm11}, $\mathcal{L}\subseteq \mathcal{L}_g$.

We now show that the inclusion is proper. Define $\sigma\in \aut$ by
\[\sigma(n)= \begin{cases}
n+x_{2m}, & \mbox{if } n\in M_{2m} \mbox{ for some } m\in \omega\\
n-x_{2m}, & \mbox{if } n\in N_{2m} \mbox{ for some } m\in \omega\\
n, & \mbox{if } n\in M_{2m+1}\cup N_{2m+1} \mbox{ for some } m\in \omega.
\end{cases}\]
We claim that $\sigma\in \mathcal L_g\setminus \mathcal L$. Fix $n\in \omega$. If $n\in M_{2m+1}\cup N_{2m+1}$ for some $m\in \omega$, then $\ks(n) = 0$. If $n\in M_{2m}\cup N_{2m}$ for some $m\in \omega$, then $\ks(n)/g(n) \leq n/n^2 = 1/n$. Thus $\ks(n)/g(n)\to 0$, and so $\sigma\in \mathcal L_g$.

On the other hand, let $n_m = y_{2m} + x_{2m}$ for every $m\in \omega$. Then \[\frac{\ks(n_m)}{n_m} = \frac{x_{2m}}{y_{2m} + x_{2m}} = \frac{x_{2m}}{2\sum_{i=1}^{2m} x_i - x_{2m}} \to 1.\] Thus $\sigma \notin \mathcal L$, which implies $\mathcal{L} \subsetneqq \mathcal{L}_g$.
\end{proof}

\begin{Ex}
\label{ex7}
There exists $g\in \gin$ for which $\mathcal L$ and $\mathcal L_g$ are incomparable.
\end{Ex}

\begin{proof}
For each $m\in \omega$, let $A_m = [2^{4m}, 2^{4m+1}) \cap \omega$, $B_m = [2^{4m+1}, 3\cdot 2^{4m}) \cap \omega$, $C_m = [2^{4m+2}, 2^{4m+2}+2^{2m+1}) \cap \omega$, and $D_m = [2^{4m+2}+2^{2m+1}, 2^{4m+2}+2^{2m+2}) \cap \omega$. Put $E_m = A_m\cup B_m$ and $F_m = C_m\cup D_m$ for all $m\in \omega$. Now define $g\in \gin$ by
\[g(n)= \begin{cases}
n^2, & \mbox{if } n\in E_m \mbox{ for some } m\in \omega\\
\sqrt n, & \mbox{if } n\in F_m \mbox{ for some } m\in \omega\\
n, & \mbox{otherwise}.
\end{cases}\]
Also define $\sigma\in \aut$ by
\[\sigma(n)= \begin{cases}
n+2^{4m}, & \mbox{if } n\in A_m \mbox{ for some } m\in \omega\\
n-2^{4m}, & \mbox{if } n\in B_m \mbox{ for some } m\in \omega\\
n, & \mbox{otherwise}.
\end{cases}\]
We show that $\sigma\in \mathcal{L}_g \setminus \mathcal{L}$. If $n\notin \bigcup_{m\in \omega} E_m$, then $\ks(n) = 0$ while if $n\in E_m$ for some $m\in \omega$, then $\ks(n)/g(n) \leq n/n^2 = 1/n$. Thus $\ks(n)/g(n)\to 0$, and so $\sigma\in \mathcal L_g$. Now for each $m\in \omega$, let $n_m = 2^{4m+1}$. Then \[\frac{\ks(n_m)}{n_m} = \frac{2^{4m}}{2^{4m+1}} \to \frac{1}{2}.\] It consequently follows that $\sigma \notin \mathcal{L}$, and hence $\sigma\in \mathcal{L}_g \setminus \mathcal{L}$. Therefore, $\mathcal{L}_g \not\subseteq \mathcal{L}$.

Next define $\tau\in \aut$ by
\[\tau(n) = \begin{cases}
n + 2^{2m+1}, & \mbox{if } n\in C_m \mbox{ for some } m\in \omega\\
n - 2^{2m+1}, & \mbox{if } n\in D_m \mbox{ for some } m\in \omega\\
n, & \mbox{otherwise}.
\end{cases}\]
We claim that $\tau\in \mathcal{L} \setminus \mathcal{L}_g$. If $n\notin \bigcup_{m\in \omega} F_m$, then $\kappa_\tau(n) = 0$. If $n\in F_m$ for some $m\in \omega$, then \[\frac{\kappa_\tau(n)}{n} \leq \frac{2^{2m+1}}{2^{4m+2}} = \frac{1}{2^{2m+1}}.\] Thus $\kappa_\tau(n)/n \to 0$, and hence $\tau\in \mathcal{L}$. Now for each $m\in \omega$, let $s_m = 2^{4m+2} + 2^{2m+1}$. Then \[ \frac{\kappa_\tau(s_m)}{g(s_m)} = \frac{2^{2m+1}}{\sqrt{2^{4m+2} + 2^{2m+1}}} \to 1.\] It follows that $\tau \notin \mathcal{L}_g$, and consequently $\tau\in \mathcal{L} \setminus \mathcal{L}_g$. Thus $\mathcal{L} \not\subseteq \mathcal{L}_g$. Therefore, $\mathcal{L}$ and $\mathcal{L}_g$ are incomparable.
\end{proof}

\begin{Prop}
\label{prop13}
Let $g_1, g_2 \in \gin$ be such that there exist $M > 0$ and $k\in \omega$ satisfying $g_1(n)/g_2(n) \leq M$ for all $n \geq k$. Then $\mathcal{L}_{g_1} \subseteq \mathcal{L}_{g_2}$.
\end{Prop}
%\begin{proof}
%Let $\sigma \in \mathcal{L}_{g_1}$. Then $\ks(n)/g_1(n)\to 0$. Now \[\frac{\ks(n)}{g_2(n)} = \frac{\ks(n)}{g_1(n)} \cdot \frac{g_1(n)}{g_2(n)} \leq M \frac{\ks(n)}{g_1(n)}\] for all $n \geq k$. It follows that $\ks(n)/g_2(n)\to 0$, and hence $\sigma \in \mathcal{L}_{g_2}$. Thus $\mathcal{L}_{g_1} \subseteq \mathcal{L}_{g_2}$.
%\end{proof}

\begin{Cor}
\label{cor5}
Let $g_1, g_2 \in \gin$ be such that there exist $M, m > 0$ and $k\in \omega$ satisfying $m \leq g_1(n)/g_2(n) \leq M$ for all $n \geq k$. Then $\mathcal{L}_{g_1} = \mathcal{L}_{g_2}$.
\end{Cor}

%We do not present any proof of this assertion as it will be clear in Section~4 that the assertion follows from the corresponding result concerning simple density ideals.

The following example shows that the condition imposed on $g_1, g_2 \in \gin$ in Corollary~\ref{cor5} is not necessary for the equality $\mathcal{L}_{g_1} = \mathcal{L}_{g_2}$. The example is taken from \cite{gd2}.

\begin{Ex}
\label{ex6}
There exist $g_1, g_2 \in \gin$ such that the sequence $(g_1(n)/g_2(n))$ is unbounded, but $\mathcal{L}_{g_1} = \mathcal{L}_{g_2}$.
\end{Ex}
\begin{proof}
Let $g_1, g_2 \in \gin$ be given by
\[g_1(n)= \begin{cases}
n, & \mbox{if } n \mbox{ is even}\\
\sqrt{n}, & \mbox{if } n \mbox{ is odd},
\end{cases} \text{\qquad and \qquad} g_2(n) = \sqrt{n}\quad \text{for all } n\in\omega.\]
Clearly, the sequence $(g_1(n)/g_2(n))$ is unbounded. Since $g_1(n)\ge g_2(n)$ for all $n\in \omega$, if $\sigma\in \mathcal{L}_{g_2}$, then $\ks(n)/g_1(n)\leq \ks(n)/g_2(n)\to 0$, which shows that $\sigma\in \mathcal{L}_{g_1}$. Thus $\mathcal{L}_{g_2} \subseteq \mathcal{L}_{g_1}$.

Conversely, let $\sigma\in \mathcal{L}_{g_1}$. Then for odd $n$, $\ks(n)/\sqrt{n} = \ks(n)/g_1(n)\to 0$. Now let $n=2m$ be even. By Lemma~\ref{lemma5}, $\ks(2m)\leq \ks(2m+1)+1$. Hence
\begin{align*}
\frac{\ks(2m)}{\sqrt{2m}} & \le \frac{\ks(2m+1)}{\sqrt{2m}}+\frac{1}{\sqrt{2m}}\\
 & = \frac{\ks(2m+1)}{\sqrt{2m+1}}\cdot \frac{\sqrt{2m+1}}{\sqrt{2m}}
+\frac{1}{\sqrt{2m}}.
\end{align*}
Since $\ks(2m+1)/\sqrt{2m+1} \to 0$ and $\sqrt{2m+1}/\sqrt{2m} \to 1$, it follows that $\ks(2m)/\sqrt{2m} \to 0$.

Combining the even and odd cases, we obtain $\ks(n)/\sqrt{n} \to 0$, and so $\sigma\in \mathcal{L}_{g_2}$. Thus $\mathcal{L}_{g_1} \subseteq \mathcal{L}_{g_2}$. Hence $\mathcal{L}_{g_1} = \mathcal{L}_{g_2}$.
\end{proof}

\begin{Rem}
\label{rem1}
Let $f\in \gin$. Suppose that $f$ is not increasing. Define a function $g: \omega \to [0, \infty)$ inductively as follows (this construction has been reproduced from \cite[Proposition 2.2]{gd2}). Choose $x_0 = \min \{f(n) : n\in \omega\}$ and $k_0 = \max \{n\in \omega : f(n) = x_0\}$. Define $g(n) = x_0$ for all $n\in k_0+1$. Next choose $x_1 = \min \{f(n) : n\notin k_0+1\}$ and $k_1 = \max \{n\in \omega : f(n) = x_1\}$. Define $g(n) = x_1$ for all $n\in (k_1+1)\setminus (k_0+1)$ and so on. Clearly, $g$ is increasing and $g(n)\to \infty$. Also $g(n) \leq f(n)$ for all $n \in \omega$. It follows that $n/f(n) \leq n/g(n)$ for all $n \in \omega$ and hence $n/g(n) \not \to 0$ as $n/f(n) \not \to 0$. Thus $g\in \gin$. In \cite[Proposition 2.2]{gd2}, authors had shown that $\mathcal{Z}_g = \mathcal{Z}_f$. Next we show that this result does not hold for simple L\'{e}vy groups (see Theorem~\ref{thm3}). Hence it is not the case that every phenomenon satisfied by simple density ideals is replicated for simple L\'{e}vy groups.
\end{Rem}

\begin{Th}
\label{thm3}
There exists $f\in \gin$ such that $f$ is not increasing but $\mathcal{L}_g \subsetneqq \mathcal{L}_f$, where $g\in \gin$ is defined as in Remark~\ref{rem1}.
\end{Th}

\begin{proof}
For each $k\geq 1$, let $M_k = [n_k, n_k+k^3]\cap \omega$ and $N_k = [n_k+k^3+1, n_k+k^3+k]\cap \omega$, where $n_1=0$ and $n_{k+1} = n_k+k^3+k+1$. Then the sequence $M_1, N_1, M_2, N_2, \dotsc$ forms a partition of $\omega$.

Define a function $f:\omega \to [0,\infty)$ by
\[f(n) = \begin{cases}
k, & \mbox{if } n=n_k \mbox{ for some } k\geq 1 \\
k^4, & \mbox{if } n\in (M_k\setminus \{n_k\}) \cup N_k \mbox{ for some } k\geq 1.
\end{cases}\]
Observe that $n_k/f(n_k) = n_k/k$ becomes unbounded as $k\to \infty$ and so $n/f(n)\not\to 0$. Also $f(n)\to \infty$ as $n\to \infty$. Thus $f\in \gin$. Clearly $f$ is not increasing.

Let $g\in \gin$ be the function defined as in Remark~\ref{rem1} associated to $f$. Then $g$ can be expressed as follows.
\[g(n) = \begin{cases}
k, & \mbox{if } n=n_k \mbox{ for some } k\geq 1 \\
k+1, & \mbox{if } n\in (M_k\setminus \{n_k\}) \cup N_k \mbox{ for some } k\geq 1.
\end{cases}\]

We claim that $\mathcal{L}_g \subsetneqq \mathcal{L}_f$. Since $(g(n)/f(n))$ is bounded, by Proposition~\ref{prop13}, we get $\mathcal{L}_g \subseteq \mathcal{L}_f$. We now define a $\sigma\in \aut$ by
\[\sigma(n) = \begin{cases}
n+k^3, & \mbox{if } n\in M_k\setminus (\{n_k\}\cup [n_k+k+1,n_k+k^3]\cap \omega) \mbox{ for some } k\geq 1\\
n, & \mbox{if }  n\in \{n_k\}\cup [n_k+k+1,n_k+k^3]\cap \omega \mbox{ for some } k\geq 1\\
n-k^3, & \mbox{if } n\in N_k \mbox{ for some } k\geq 1.
\end{cases}\]
Thus we can say that for each $k\geq 1$, $\sigma(\{n_k\}\cup [n_k+k+1,n_k+k^3]\cap \omega) = N_k$, $\sigma(N_k) = M_k\setminus (\{n_k\}\cup [n_k+k+1,n_k+k^3]\cap \omega)$ and $\sigma(n) = n$ on $\{n_k\}\cup [n_k+k+1,n_k+k^3]\cap \omega$. We first show that $\sigma\in \mathcal{L}_f$. Note that for  $n\in (M_k\setminus\{n_k\}) \cup N_k$, $k\geq 1$, we have
\[\frac{\ks(n)}{f(n)} \leq \frac{k}{k^4} = \frac{1}{k^3} \to 0.\] It then follows that $\ks(n)/f(n)\to 0$ as $n\to \infty$ along $\bigcup_{k\geq 1}((M_k\setminus\{n_k\}) \cup N_k)$.

If $n= n_k$, $k\geq 1$, then $\ks(n) = 0$. This gives us $\ks(n)/f(n) \to 0$ as $n\to \infty$ along $\{n_k : k\geq 1\}$. Thus $\ks(n)/f(n) \to 0$, and hence $\sigma\in \mathcal{L}_f$.

We finally show that $\sigma \not\in \mathcal{L}_g$. Let $n = n_k + 1 + \lfloor k^3/2\rfloor$, $k\geq 1$. Then \[\frac{\ks(n)}{g(n)} = \frac{k}{k+1} \to 1\] and hence  $\ks(n)/g(n)\to 1$ as $n\to \infty$ along $\{n_k + 1 + \lfloor k^3/2\rfloor : k\geq 1\}$. Thus $\sigma\notin \mathcal{L}_g$. Therefore, $\mathcal{L}_g \subsetneqq \mathcal{L}_f$.
\end{proof}

\begin{Th}
\label{thm4}
Let $g_1,g_2\in \gin$ be such that $g_2(n)/g_1(n)\to \infty$. Further assume that $g_2$ is eventually increasing and that there exist $c>0$ and $d\in \omega$ such that $g_1(n+1)-g_1(n)\leq c$ for every $n\geq d$. Then $\mathcal{L}_{g_1} \subsetneqq \mathcal{L}_{g_2}$.
\end{Th}

\begin{proof}
Since $g_1(n)/g_2(n)\to 0$, by Proposition~\ref{prop13}, we have $\mathcal{L}_{g_1} \subseteq \mathcal{L}_{g_2}$. We have to show that the inclusion is proper. Since $g_2(n)/g_1(n)\to \infty$ and $g_1(n)\to \infty$, we may choose a strictly increasing sequence $(n_k)$ in $\omega$ such that $n_k\geq d$ for all $k\in \omega$, and the intervals $I_k\cup J_k$ are pairwise disjoint, where for each $k\in \omega$, $I_k = [n_k, n_k+\ell_k) \cap \omega$, $J_k = [n_k+\ell_k, n_k+2\ell_k) \cap \omega$, and $\ell_k = \lfloor g_1(n_k)\rfloor$.

Define $\sigma\in \aut$ by
\[\sigma(n) = \begin{cases}
n+\ell_k, & \mbox{if } n\in I_k \mbox{ for some } k\in \omega\\
n-\ell_k, & \mbox{if } n\in J_k \mbox{ for some } k\in \omega\\
n, & \mbox{otherwise}.
\end{cases}\]
We first show that $\sigma\in \mathcal{L}_{g_2}$. If $n\notin \bigcup_{k\in \omega} (I_k\cup J_k)$, then $\ks(n) = 0$. Let $n\in \omega$ be large enough. If $n\in I_k\cup J_k$ for some $k\in \omega$, then $\ks(n)\leq \ell_k$. Since $g_2$ is eventually increasing and $n\geq n_k$ on $I_k\cup J_k$, we obtain $g_2(n)\geq g_2(n_k)$. Thus \[\frac{\ks(n)}{g_2(n)} \leq \frac{\ell_k}{g_2(n_k)} \leq \frac{g_1(n_k)}{g_2(n_k)} \to 0.\] Therefore, $\sigma\in \mathcal{L}_{g_2}$.

We now show that $\sigma \notin \mathcal{L}_{g_1}$. Put $m_k = n_k+\ell_k$ for every $k\in \omega$. Then $\ks(m_k) = \ell_k$ for every $k\in \omega$. Since $g_1(n+1)-g_1(n)\leq c$ for all $n\geq d$, we have $g_1(m_k)\leq g_1(n_k) + c\ell_k$ for every $k\in \omega$. Since $\ell_k = \lfloor g_1(n_k)\rfloor$, it follows that \[g_1(m_k)\leq g_1(n_k) + c g_1(n_k) = (1+c)g_1(n_k)\] for every $k\in \omega$. Also, since $g_1(n_k)\to \infty$, we have $\ell_k = \lfloor g_1(n_k)\rfloor \geq g_1(n_k)/2$ for all sufficiently large $k\in \omega$. Therefore, \[\frac{\ks(m_k)}{g_1(m_k)} = \frac{\ell_k}{g_1(m_k)} \geq \frac{g_1(n_k)/2}{(1+c)g_1(n_k)} = \frac{1}{2(1+c)}\] for all sufficiently large $k\in \omega$. Hence $\ks(n)/g_1(n) \not\to 0$, and so $\sigma \notin \mathcal{L}_{g_1}$. Thus $\mathcal{L}_{g_1} \subsetneqq \mathcal{L}_{g_2}$.
\end{proof}

Example~\ref{ex7} establishes the existence of a single function $g\in \gin$ such that $\mathcal{L}_g$ and $\mathcal{L}$ are incomparable. We now strengthen this result by constructing a continuum-sized family of such functions, with the additional property that the corresponding simple L\'{e}vy groups are pairwise incomparable.

%Though the following result follows from \cite[Theorem~2.7]{gd2}, in view of the interrelations established between simple density ideals and simple L\'{e}vy groups in Section~4, but one can always ask if a direct proof can be given and we proceed exactly in that direction and present a direct proof without bypassing the notion of ideals. Further it should be noted that Theorem 2.7 \cite{gd2}, which in some sense is the main result of that article also follows from the theorem established below. So in some sense it is like looking at some object from two different perspectives.  As a result, one must understand that the essential tool is same in both results, namely, almost disjoint ness. Recall that a family $\mathcal{A}\subseteq [\omega]^\omega$ is said to be almost disjoint if $|A\cap B|<\omega$ for any two distinct $A,B\in \mathcal{A}$ (where $[\omega]^\omega$ denotes the family of all infinite subsets of $\omega$).

\begin{Th}
\label{thm8}
There exists a family $G\subseteq \gin$ of cardinality $\mathfrak{c}$ such that, for every $g\in G$, $\mathcal{L}_g$ and $\mathcal{L}$ are incomparable, and for every distinct $f,g\in G$, $\mathcal{L}_f$ and $\mathcal{L}_g$ are incomparable.
\end{Th}

\begin{proof}
For each $k\in \omega$, let $m_k = 2^{2^{k+2}}$ and $B_k = [m_k,5m_k) \cap \omega$. Then $(B_k)$ is a pairwise disjoint family of intervals, since $5m_k<m_{k+1}$ for every $k\in \omega$.

Let $\mathcal{A}\subseteq [\omega]^\omega$ be an almost disjoint family of cardinality $\mathfrak{c}$. For each $A\in \mathcal{A}$, define a function $g_A:\omega\to [0,\infty)$ by
\[g_A(n)= \begin{cases}
k, & \mbox{if } n\in B_k \mbox{ and } k\in A\\
m_k^2, & \mbox{if } n\in B_k \mbox{ and } k\notin A\\
n, & \mbox{if } n\notin \bigcup_{k\in \omega} B_k.
\end{cases}\]
Put $G = \{g_A : A\in \mathcal{A}\}$. We first show that $G\subseteq \gin$. Fix $A\in \mathcal{A}$. Clearly, $g_A(n)\to \infty$ as $n\to \infty$. Moreover, since $A$ is infinite and for every $k\in A$ and every $n\in B_k$, $g_A(n)=k$, it follows that \[\frac{m_k}{g_A(m_k)} = \frac{m_k}{k}\to \infty \text{ as } k\to \infty \text{ along } A.\]
Hence $n/g_A(n)\not\to 0$. Thus $g_A\in \gin$, and therefore $G\subseteq \gin$. Again as the map $A\mapsto g_A$ is bijective, we further have $|G| = |\mathcal{A}| = \mathfrak{c}$. %The map $A\mapsto g_A$ is injective. Indeed, if $A,B\in \mathcal{A}$ are distinct, then there exists $k\in \omega$ such that $k\in A\setminus B$. Hence $g_A(m_k) = k \neq m_k^2 = g_B(m_k)$, and therefore $g_A \neq g_B$. By definition, the map $A\mapsto g_A$ is surjective.

Next, for every infinite $E\subseteq \omega$, define $\sigma_E, \tau_E\in \aut$ by \[\sigma_E(n) = \begin{cases}
n+k, & \mbox{if } n\in [m_k,m_k+k) \cap \omega \mbox{ and } k\in E\\
n-k, & \mbox{if } n\in [m_k+k,m_k+2k) \cap \omega \mbox{ and } k\in E\\
n, & \mbox{otherwise},
\end{cases}\]
and
\[\tau_E(n) = \begin{cases}
n+m_k, & \mbox{if } n\in [m_k,2m_k) \cap \omega \mbox{ and } k\in E\\
n-m_k, & \mbox{if } n\in [2m_k,3m_k) \cap \omega \mbox{ and } k\in E\\
n, & \mbox{otherwise}.
\end{cases}\]

We now record two simple facts Claim~1 and Claim~2.

\medskip
\noindent \textbf{Claim 1.} For every infinite $E\subseteq \omega$, one has $\sigma_E\in \mathcal{L}$.

\smallskip
Indeed, if $n\notin [m_k,m_k+2k) \cap \omega$ for all $k\in E$, then $\kappa_{\sigma_E}(n) = 0$ while  $n\in [m_k,m_k+2k) \cap \omega$ for some $k\in E$ implies that $\kappa_{\sigma_E}(n) \leq k$, and so \[\frac{\kappa_{\sigma_E}(n)}{n} \leq \frac{k}{m_k}.\] Thus $\kappa_{\sigma_E}(n)/n \to 0$, and hence $\sigma_E\in \mathcal{L}$.

\medskip
\noindent
\textbf{Claim 2.} Let $A\in \mathcal{A}$ and let $E\subseteq \omega$ be infinite.
\begin{enumerate}
\item[(a)] If $E\subseteq \omega\setminus A$, then $\sigma_E\in \mathcal{L}_{g_A}$.
\item[(b)] If $E\subseteq A$, then $\sigma_E\notin \mathcal{L}_{g_A}$.
\end{enumerate}

\smallskip
To prove (a), assume that $E\subseteq \omega\setminus A$. If $n\in [m_k,m_k+2k) \cap \omega$ for some $k\in E$, then $g_A(n) = m_k^2$ and $\kappa_{\sigma_E}(n)\leq k$, and so \[\frac{\kappa_{\sigma_E}(n)}{g_A(n)} \leq \frac{k}{m_k^2}.\] Note that outside these intervals the quotient is zero. Thus $\kappa_{\sigma_E}(n)/g_A(n) \to 0$, i.e. $\sigma_E\in \mathcal{L}_{g_A}$.

To prove (b), assume that $E\subseteq A$. For each $k\in E$, let $n_k = m_k+k$. Then for each $k\in E$, $\kappa_{\sigma_E}(n_k) = k$ and $g_A(n_k) =k$, and therefore $\kappa_{\sigma_E}(n_k)/g_A(n_k) = 1$. Thus $\sigma_E\notin \mathcal{L}_{g_A}$.

\medskip
\noindent
\textbf{Claim 3.} Let $A\in \mathcal{A}$ and let $E\subseteq \omega\setminus A$ be infinite. Then $\tau_E\in \mathcal{L}_{g_A}\setminus \mathcal{L}$.

\smallskip
For each $k\in \omega$, let $n_k = 2m_k$. Then for each $k\in E$, $\kappa_{\tau_E}(n_k) = m_k$, and so \[\frac{\kappa_{\tau_E}(n_k)}{n_k} = \frac{m_k}{2m_k} = \frac{1}{2}.\] Hence $\tau_E\notin \mathcal{L}$.

On the other hand, if $n\in [m_k,3m_k) \cap \omega$ for some $k\in E$, then $\kappa_{\tau_E}(n) \leq m_k$, and since $k\notin A$, we have $g_A(n) = m_k^2$. It follows that \[\frac{\kappa_{\tau_E}(n)}{g_A(n)} \leq \frac{m_k}{m_k^2} = \frac{1}{m_k}.\] Again observe that outside these intervals the quotient is zero. Therefore, $\kappa_{\tau_E}(n)/g_A(n) \to 0$, i.e. $\tau_E\in \mathcal{L}_{g_A}$, proving the claim.

We now verify the required incomparability statements. Fix $A\in \mathcal{A}$. Since $A$ is infinite, Claim~1 implies $\sigma_A\in \mathcal{L}$, while Claim~2(b) yields $\sigma_A\notin \mathcal{L}_{g_A}$. Hence $\mathcal{L} \nsubseteq \mathcal{L}_{g_A}$. Also since $\omega\setminus A$ is infinite, Claim~3 applied to $E = \omega\setminus A$ gives $\tau_{\omega\setminus A} \in \mathcal{L}_{g_A} \setminus \mathcal{L}$, and so $\mathcal{L}_{g_A} \nsubseteq \mathcal{L}$. Thus $\mathcal{L}_{g_A}$ and $\mathcal{L}$ are incomparable.

Finally, let $A,B\in \mathcal{A}$ be distinct. Since $\mathcal{A}$ is almost disjoint, the sets $A\setminus B$ and $B\setminus A$ are both infinite. Put $E = A\setminus B$. Then $E\subseteq A$ and $E\subseteq \omega\setminus B$. By Claim~2(b), $\sigma_E\notin \mathcal{L}_{g_A}$, while by Claim~2(a), $\sigma_E\in \mathcal{L}_{g_B}$. Hence $\mathcal{L}_{g_B} \nsubseteq \mathcal{L}_{g_A}$. Similarly, taking $F = B\setminus A$, we obtain $\sigma_F\in \mathcal{L}_{g_A}\setminus \mathcal{L}_{g_B}$, and therefore $\mathcal{L}_{g_A} \nsubseteq \mathcal{L}_{g_B}$. Thus $\mathcal{L}_{g_A}$ and $\mathcal{L}_{g_B}$ are incomparable.

Consequently, the family $G$ has all the required properties.
\end{proof}

\subsection{Union and intersection of simple L\'{e}vy groups}
It is natural to ask what is exactly the intersection and union of all possible simple L\'{e}vy groups and the following two results answer the same.

\begin{Prop}
\label{propn2}
$\bigcap_{g\in \gin} \mathcal{L}_g = \{\sigma\in \aut : \kappa_\sigma \text{ is bounded}\}$.
\end{Prop}

\begin{proof}
Put $\mathcal{B} = \{\sigma\in \aut : \kappa_\sigma \text{ is bounded}\}$. Let $\sigma\in \mathcal{B}$. Then there exists $C>0$ such that $\ks(n)\leq C$ for all $n\in \omega$. If $g\in \gin$, then $g(n)\to \infty$, and hence $\ks(n)/g(n)\leq C/g(n)\to 0$. Thus $\sigma\in \bigcap_{g\in\gin} \mathcal{L}_g$, and therefore $\mathcal{B} \subseteq \bigcap_{g\in\gin} \mathcal{L}_g$.

Conversely, let $\sigma\in \bigcap_{g\in\gin} \mathcal{L}_g$. Suppose that $\ks$ is unbounded. For each $n\in \omega$, choose $f(n) = \max\{\ks(k) : k\in n+1\}$. Since $\ks$ is unbounded, $f(n)\to \infty$. Moreover, there exists a strictly increasing sequence $(n_j)$ in $\omega$ such that $\ks(n_j) = f(n_j)\to \infty$. Now define $g\in \gin$ by $g(n) = \sqrt{f(n)}$ for all $n\in \omega$. Then \[\frac{\ks(n_j)}{g(n_j)} = \frac{f(n_j)}{\sqrt{f(n_j)}} = \sqrt{f(n_j)} \to \infty.\] It follows that $\sigma \notin \mathcal{L}_g$, which is a contradiction. Thus $\ks$ is bounded, i.e. $\sigma \in \mathcal{B}$, and hence $\bigcap_{g\in\gin} \mathcal{L}_g \subseteq \mathcal{B}$. Therefore, $\bigcap_{g\in\gin} \mathcal{L}_g = \mathcal{B}$.
\end{proof}

\begin{Prop}
\label{propn3}
$\bigcup_{g\in \gin} \mathcal{L}_g = \left\{\sigma\in \aut : \liminf_n \frac{\ks(n)}{n} = 0\right\}$. Moreover, $\bigcup_{g\in \gin} \mathcal{L}_g \subsetneqq \aut$.
\end{Prop}

\begin{proof}
Put $\mathcal{B} = \left\{\sigma\in \aut : \liminf_n \frac{\ks(n)}{n} = 0\right\}$. Let $\sigma\in \bigcup_{g\in \gin} \mathcal{L}_g$. Then $\sigma\in \mathcal{L}_g$ for some $g\in \gin$, and hence $\ks(n)/g(n)\to 0$. Since $n/g(n)\not\to 0$, there exist $\varepsilon >0$ and an increasing sequence $(n_j)$ in $\omega$ such that $n_j/g(n_j)\geq \varepsilon$ for all $j\in \omega$. Therefore, \[\frac{\ks(n_j)}{n_j} = \frac{\ks(n_j)}{g(n_j)} \cdot \frac{g(n_j)}{n_j} \leq \frac{\ks(n_j)}{\varepsilon g(n_j)} \to 0.\] Consequently, $\liminf_n \frac{\ks(n)}{n} = 0$, i.e. $\sigma\in \mathcal{B}$. Hence $\bigcup_{g\in \gin} \mathcal{L}_g \subseteq \mathcal{B}$.

Conversely, let $\sigma\in \mathcal{B}$. Then $\liminf_n \frac{\ks(n)}{n} = 0$. Choose an increasing sequence $(n_j)$ in $\omega$ such that $n_j\geq j^2$ and $\ks(n_j)/n_j\to 0$. Define $g\in \gin$ by
\[g(n) = \begin{cases}
\sqrt{n_j(\ks(n_j)+1)}, & \mbox{if } n = n_j \mbox{ for some } j\in \omega\\
(n+1) (\ks(n)+1), & \mbox{otherwise}.
\end{cases}\]
If $n\notin \{n_j : j\in \omega\}$, then \[\frac{\ks(n)}{g(n)} = \frac{\ks(n)}
{(n+1)(\ks(n)+1)} \leq \frac{1}{n+1}.\] On the other hand, \[\frac{\ks(n_j)}{g(n_j)} = \frac{\ks(n_j)}{\sqrt{n_j(\ks(n_j)+1)}} \leq \sqrt{\frac{\ks(n_j)}{n_j}}.\] Thus $\ks(n)/g(n)\to 0$, i.e. $\sigma\in \mathcal{L}_g$. Hence $\mathcal{B} \subseteq \bigcup_{g\in \gin} \mathcal{L}_g$. Therefore, $\bigcup_{g\in \gin} \mathcal{L}_g = \mathcal{B}$.

It remains to show that $\bigcup_{g\in \gin} \mathcal{L}_g \subsetneqq \aut$. Put $B = \{4k+3 : k\in \omega\}$ and $A = \omega \setminus B$, and let $(a_k)$ be the increasing enumeration of $A$. Define $\sigma\in \aut$ by
\[\sigma(n) = \begin{cases}
4k+3, & \mbox{if } n = a_k \mbox{ for some } k\in \omega\\
a_k, & \mbox{if } n = 4k+3 \mbox{ for some } k\in \omega.
\end{cases}\]
Since $a_k< 4k+3$ for all $k\in \omega$, we have $\ks(n) = |A\cap n| - |B\cap n| = n - 2|B\cap n|$ for all $n\in \omega$. Moreover, $|B\cap n| = \lfloor n/4\rfloor$, and hence $\ks(n) = n - 2\lfloor n/4\rfloor \geq n/2$ for all $n\in \omega$. Therefore, $\liminf_n \frac{\ks(n)}{n} \geq 1/2$. By the first part of the proof, $\sigma \notin \bigcup_{g\in \gin} \mathcal{L}_g$. Thus $\bigcup_{g\in \gin} \mathcal{L}_g \subsetneqq \aut$.
\end{proof}

\vspace{1 cm}

\noindent\textbf{Statements and Declarations:}\\

\noindent The authors declare that no data was used for the preparation of this manuscript.\\

\noindent The authors declare that no funds, grants, or other support were received during the preparation of this manuscript.\\

\noindent The authors have no relevant financial or non-financial interests to disclose.\\

\noindent All authors contributed to the study conception, material preparation and analysis. The first draft of the manuscript was written by Nur Alam  and all authors commented on previous versions of the manuscript. All authors read and approved the final manuscript.
\end{document}